\documentclass[conference,10pt]{IEEEtran}

\usepackage{amsmath,amssymb,amsfonts,amsthm}
\usepackage{bm}
\usepackage{cite}
\usepackage{microtype}
\newtheorem{theorem}{Theorem}
\newtheorem{lemma}[theorem]{Lemma}
\newtheorem{proposition}[theorem]{Proposition}
\newtheorem{corollary}[theorem]{Corollary}
\newtheorem{definition}[theorem]{Definition}
\newtheorem{remark}[theorem]{Remark}
\newtheorem{observation}[theorem]{Observation}

\newcommand{\F}{\mathbb{F}}

\newcommand{\wt}{\operatorname{wt}}
\newcommand{\supp}{\operatorname{supp}}
\newcommand{\Aw}{\mathcal{A}}
\newcommand{\C}{\mathcal{C}}
\newcommand{\PG}{\operatorname{PG}}

\begin{document}
	
	\title{A Complete Characterization of \(q\)-Ary \(2\)-Designs in Doubly-Extended Reed-Solomon Codes}
	
	\author{%
		\IEEEauthorblockN{Xinghao Wu and Junling Zhou}
		\IEEEauthorblockA{School of Mathematics and Statistics, Beijing Jiaotong University\\
			Beijing 100044, China\\
			Emails: 23111512@bjtu.edu.cn, jlzhou@bjtu.edu.cn}
	}
	
	\maketitle
	
	\begin{abstract}
Let \({\cal C}_k\) be the doubly-extended Reed--Solomon code with parameters \([q+1, k, q-k+2]_q\). We give a complete characterization of when its codewords of a fixed weight form \(q\)-ary 2-designs. The main ingredient is an exact criterion showing that the number of \(r\)-subsets of \(\F_q^*\) with a prescribed product is independent of that product if and only if \(\gcd(r, q-1)=1\). Applying this criterion to the minimum-weight codewords of \({\cal C}_k\), we prove that they form a \(q\)-ary 2-design if and only if \(\gcd(k-1, q-1)=1\). We further establish a propagation result for MDS codes showing that the \(q\)-ary 2-design property at minimum weight implies the same property at every nonzero weight. Consequently, \({\cal C}_k\) holds \(q\)-ary 2-designs at all nonzero weights if and only if \(\gcd(k-1, q-1)=1\).
	\end{abstract}
	
	
	\section{Introduction}
	The interaction between coding theory and design theory has been extensively studied for several decades. On the one hand, incidence matrices of combinatorial designs provide a natural source of linear codes. On the other hand, the supports of codewords of a fixed weight in a linear code may yield combinatorial designs. A fundamental result in the latter direction is the Assmus--Mattson theorem \cite{AM1969}. A large body of subsequent work has developed constructions, characterizations, and infinite families of designs arising from linear and nonlinear codes; see, for example, \cite{Key,Ding,BCH-4-design,Tonchev1,Tonchev2,Ding-DS,Awada,Cunsheng2020Infinite,Heng,Yan} and the references therein. Thus, the theory of ordinary designs arising from codes is by now well established.
	
	A vector analogue of ordinary designs, now called a $q$-ary $t$-design, was introduced in the early work of Delsarte \cite{dsd's',Delsarte1973}. In a $q$-ary design, blocks are vectors over $\F_q$ rather than subsets of coordinate positions, so both the locations and the values of their nonzero entries are retained. Consequently, every $q$-ary $t$-design induces an ordinary $t$-design through its supports, whereas the converse does not hold in general. Early results of Assmus, Goethals, and Mattson \cite{punctured-t} and Goethals and van Tilborg \cite{Tilborg1975,UP} revealed connections between $q$-ary designs, regularity, and linear codes. In contrast with the extensive literature on ordinary designs supported by codes, however, comparatively few results have been obtained on $q$-ary $t$-designs arising directly from codes over nonbinary fields.
	
	In our recent work \cite{WuZhou2026}, we initiated a systematic investigation of linear codes holding $q$-ary $t$-designs. General criteria based on the parameters and regularity of a code were developed, together with an automorphism-group approach for cases beyond the scope of these criteria. In particular, for the doubly-extended Reed--Solomon (DRS) code $\C_k$ with parameters $[q+1,k,q-k+2]_q$, it was proved that its minimum-weight codewords form a $q$-ary $2$-design whenever
	\[
	\gcd(k-1,q-1)=1.
	\]
	This naturally raises two questions: whether this coprimality condition is also necessary, and whether the same condition determines the $q$-ary $2$-design property at all nonzero weights.
	
	In this paper, we answer both questions completely. First, we establish an exact criterion for when the number of fixed-size subsets of the cyclic group $\F_q^*$ with a prescribed product is independent of the prescribed product. This criterion shows that the above coprimality condition is not only sufficient but also necessary for the minimum-weight codewords of a DRS code to form a $q$-ary $2$-design. Second, we prove a propagation result for MDS codes: if the minimum-weight codewords form a $q$-ary $2$-design, then so do the codewords of every nonzero weight. Consequently, $\gcd(k-1,q-1)=1$
	is equivalent both to the $q$-ary $2$-design property at minimum weight and to the corresponding property at every nonzero weight of $\C_k$.
	
	The rest of the paper is organized as follows. Section~II recalls the necessary definitions and basic properties of DRS codes. Section~III establishes the propagation result for MDS codes. Section~IV proves the subset-product criterion and derives the complete characterization of DRS codes. Section~V concludes the paper.

	\section{Preliminaries}
	
	Let $q$ be a prime power, and let $\F_q$ denote the finite field with $q$ elements. For a vector $\bm{x}=(x_1,x_2,\ldots,x_n)\in\F_q^n$, its \emph{Hamming weight} and \emph{support} are defined by $\wt(\bm{x})=|\{i:x_i\neq 0\}|$, $\supp(\bm{x})=\{i:x_i\neq 0\}$, respectively. The Hamming distance between $\bm{x},\bm{y}\in\F_q^n$ is $d(\bm{x},\bm{y})=\wt(\bm{x}-\bm{y})$.
	
	An $[n,k,d]_q$ linear code $\C$ is a $k$-dimensional subspace of $\F_q^n$ with minimum distance $d$. For each integer $w$, let
	\[
	\Aw_w(\C)=\{\bm{c}\in\C:\wt(\bm{c})=w\},\quad A_w=|\Aw_w(\C)|.
	\]
	The code $\C$ is called \emph{maximum distance separable} (MDS) if $d=n-k+1$.
	
	A $t$-$(n,w,\lambda)$ \emph{design} is a pair $(\mathcal{P},\mathcal{B})$, where $\mathcal{P}$ is an $n$-element point set and $\mathcal{B}$ is a family of $w$-subsets of $\mathcal{P}$, called \emph{blocks}, such that every $t$-subset of $\mathcal{P}$ is contained in exactly $\lambda$ blocks. The parameters $t$ and $\lambda$ are called the \emph{strength} and \emph{index}, respectively. When the point set is clear from the context, we identify the design with its block set $\mathcal{B}$.
	
	Denote by $\binom{\mathcal{P}}{w}$ the set of all $w$-subsets of $\mathcal{P}$. Then $\binom{\mathcal{P}}{w}$ forms the complete $w$-$(n,w,1)$ design. Moreover, every $t$-$(n,w,\lambda)$ design is also an $i$-$(n,w,\lambda_i)$ design for each $0\leq i\leq t$.
	To distinguish these designs from their $q$-ary analogues introduced below, they will also be referred to as \emph{ordinary} or \emph{classical} designs.
	
	The notion of a classical design can be extended from subsets to vectors over finite fields, leading to $q$-ary designs; see, for example, \cite{punctured-t,dsd's',Delsarte1973,WuZhou2026}. For $\bm{a},\bm{b}\in\F_q^n$, we say that $\bm{a}$ \emph{covers} $\bm{b}$ if $a_i=b_i$ whenever $b_i\neq 0 $.
		
	\begin{definition}\label{def:qdesign}
		Let $\Aw$ be a nonempty set of vectors of constant weight $w$ in $\F_q^n$, whose elements are called \emph{blocks}. Then $\Aw$ is a \emph{$q$-ary $t$-$(n,w,\lambda)$ design} if every vector of weight $t$ is covered by exactly $\lambda$ blocks of $\Aw$. Equivalently, for every $\bm{x}\in\F_q^n$ with $\wt(\bm{x})=t$, there are exactly $\lambda$ vectors $\bm{a}\in\Aw$ such that $d(\bm{a},\bm{x})=w-t$.
	\end{definition}
	
	We now describe the two types of designs arising from codewords of a fixed weight. Index the coordinates of $\C$ by an $n$-element set $\mathcal{P}$. For each $w$ with $A_w\neq 0$, define
	\[
	\mathcal{B}_w(\C)=\{\supp(\bm{c}):\bm{c}\in\Aw_w(\C)\}.
	\]
	Thus, $\Aw_w(\C)$ consists of the codewords of weight $w$ themselves, whereas $\mathcal{B}_w(\C)$ records only their supports. If $\Aw_w(\C)$ is a $q$-ary $t$-$(n,w,\lambda)$ design, we say that $\C$ \emph{holds} a $q$-ary $t$-design at weight $w$. On the other hand, if $ (\mathcal{P},\mathcal{B}_w(\C))$ is a classical $t$-$(n,w,\lambda)$ design, then $\C$ is said to \emph{support a $t$-design at weight $w$}.

	Assume that $\Aw_w(\C)$ forms a $q$-ary $t$-$(n,w,\lambda)$ design. Double counting the pairs $(\bm{x},\bm{c})$ satisfying $\wt(\bm{x})=t,  \bm{c}\in\Aw_w(\C), \bm{c}\text{ covers }\bm{x}$  gives
	\begin{equation}\label{eq:lambda-double-count}
		\lambda (q-1)^t\binom{n}{t}=A_w\binom{w}{t}.
	\end{equation}
	Consequently, $\lambda$ is determined by the weight distribution of $\C$. Moreover,  for every integer $i$ with $0\le i\le t$, $\Aw_w(\C)$ is also a $q$-ary $i$-$(n,w,\lambda_{i})$ design, where
	\begin{align}
		\lambda_{i}=\lambda(q-1)^{t-i}\binom{n-i}{t-i}\bigg/\binom{w-i}{t-i}.
	\end{align}
	Therefore, a necessary condition for the existence of such a design is that all $\lambda_{i}$ are integers.

	We next recall DRS codes. Let $\alpha$ be a primitive element of $\F_q$ and let
	\[
	\mathcal{P}_k[x]=\{f\in\F_q[x]:\deg f\le k-1\},\qquad 1\le k\le q.
	\]
	For $f\in \mathcal{P}_{k}[x]$, let $f(\infty)$ denote the coefficient of $x^{k-1}$.
	The DRS code is
	\begin{equation}\label{eq:DRS}
		\C_k=\bigl\{(f(0),f(1),f(\alpha),\ldots,f(\alpha^{q-2}),f(\infty)):
		f\in \mathcal{P}_k[x]\bigr\}.
	\end{equation}
	It is the projective Reed--Solomon code obtained by evaluating at the
	$q+1$ points of $\PG(1,q)$ and has parameters
	\begin{equation}\label{eq:DRSpar}
		[q+1,k,d=q-k+2]_q.
	\end{equation}
	In particular, when $k=1$, $\mathcal{P}_1[x]$ consists of constant
	polynomials, and hence
	\[
	\C_1
	=
	\{(a,a,\ldots,a):a\in\F_q\},
	\]
	which is the $[q+1,1,q+1]_q$ repetition code.
	
	The family of DRS codes is closed under duality. More precisely, it
	follows from \cite[Corollary~4.11]{Gimenez2023PRS} that
	\begin{equation}
		\label{eq:DRS-dual}
		\C_k^\perp=\C_{q+1-k},
		\qquad 1\le k\le q.
	\end{equation}
	Consequently, the boundary case $k=q$ gives $\C_q^\perp=\C_1$,
	so the dual code is again contained in the same DRS family.
	
	In the remainder of this paper, unless otherwise stated, we focus on
	the nontrivial range $2\le k\le q$,
	for which the minimum distance of $\C_k$ is $d=q-k+2$.

	For an MDS code, every $d$-subset of the coordinate set occurs as the support of a minimum-weight codeword, and each support corresponds to exactly $q-1$ nonzero scalar multiples \cite[Ch.~7]{HuffmanPless}. Hence
	\begin{equation}\label{eq:Ad-DRS}
		|\Aw_d(\C_k)|=(q-1)\binom{q+1}{d}.
	\end{equation}
	If these words form a $q$-ary $2$-design, substituting
	\eqref{eq:Ad-DRS} into \eqref{eq:lambda-double-count} already forces
	\begin{equation}\label{eq:lambda-forced}
		\lambda=\frac{1}{q-1}\binom{q-1}{k-1}.
	\end{equation}	
	Since every $d$-subset occurs as a support,
	$\mathcal{B}_d(\C_k)$ is the complete $d$-$(q+1,d,1)$ design.
	Generally, it is a	$t$-$(q+1,d,\binom{q+1-t}{d-t})$ design for every $1\le t\le d$.  However, $\Aw_d(\C_k)$ is not necessarily a $q$-ary $2$-design. We aim to characterize the $q$-ary design properties in the DRS codes in this paper.

	The characterization in Section~IV will reduce the covering problem for minimum-weight DRS codewords to a subset-product problem over $\F_q^{*}=\F_{q}\setminus\{0\}$. We therefore record the basic facts about multiplicative	characters needed for that argument.
	
	A \emph{multiplicative character} of $\F_q^{*}$ is a group homomorphism
	\[
	\chi:\F_q^{*}\longrightarrow\mathbb{C}^{*}.
	\]
	The trivial multiplicative character, denoted by $\chi_0$, is defined by $\chi_0(x)=1$ for all $x\in\F_q^{*}$. The order of a multiplicative character $\chi$ is the least positive integer $d$ such that $\chi^d=\chi_0$.
	
	Since $\F_q^{*}$ is cyclic of order $q-1$, a multiplicative character of order $d$ exists if and only if $d\mid(q-1)$. If $\chi$ has order $d$, then $\chi(\F_q^{*})=U_d$, where $U_d$ denotes the group of complex $d$-th roots of unity. By the first isomorphism theorem, $|\ker(\chi)|=(q-1)/d$,
	so every element of $U_d$ occurs exactly $(q-1)/d$ times among the values $\chi(x)$ as $x$ ranges over $\F_q^{*}$.
	
	The standard orthogonality relation for nontrivial multiplicative characters will be used below.
	\begin{lemma}\cite[Ch.~5]{LidlNiederreiter}\label{mul-char-sum}
		Let $\chi$ be a nontrivial multiplicative character of $\F_q^{*}$. Then $\sum_{x\in\F_q^{*}}\chi(x)=0$.
	\end{lemma}
	
	\section{A Propagation Result for MDS Codes}
	
	To extend the $q$-ary $2$-design property from the minimum-weight to all nonzero weights in an MDS code, we first	derive a local counting formula for $q$-ary designs.
	
	\begin{lemma}\cite[Lemma 3.10]{WuZhou2026}	\label{lem:local-pattern}
		Let $\Aw\subseteq\F_q^n$ be a $q$-ary
		$t$-$(n,w,\lambda)$ design.
		Let $A,B,C\subseteq\{1,2,\dots,n\}$ be pairwise disjoint subsets with $|A|=a,\ |B|=b,\ |C|=c$, where $ a+b+c\le t$.
		For each $i\in A\cup B$, fix $\alpha_i\in\F_q^*$. Then the number of vectors $\bm{v}\in\Aw$ satisfying
		\[
		v_i=\alpha_i \ (i\in A),\ v_i\in\F_q^\ast\setminus\{\alpha_i\}\ (i\in B),\ v_i=0\ (i\in C)
		\]
		is independent of the particular choice of the $\alpha_i$s and depends only on $a,b,c$. Denote this number by $\lambda_{a,b,c}$. More precisely,
		\begin{align*}
			\lambda_{a,b,c}=\lambda(q-2)^b (q-1)^{t-a-b}\frac{\binom{n-a-b-c}{w-a-b}}{\binom{n-t}{w-t}}.
		\end{align*}
	\end{lemma}
	We next convert these local coordinate counts into regularity of the distance distribution.
	For an $[n,k,d]_{q}$ linear code $\C$ and a vector $\bm{x}\in\F_q^n$ define the outer distribution
	\[
	B_{\bm{x},j}=|\{\bm{c}\in\C:d(\bm{x},\bm{c})=j\}|,
	\]
	and let $B_{\bm{x},j}(w)$ denote the contribution from codewords of weight
	$w$.  We need the following elementary local consequence of the design
	property.
	
	\begin{lemma}\label{lem:local-distance}
		Suppose $\Aw_w(\C)$ is a $q$-ary $2$-design. If
		$\bm{x},\bm{y}\in\F_q^n$ satisfy $\wt(\bm{x})=\wt(\bm{y})\in\{1,2\}$, then
		\[
		B_{\bm{x},j}(w)=B_{\bm{y},j}(w)\qquad(0\le j\le n).
		\]
	\end{lemma}
	\begin{proof}
		Let $\bm{x}\in\F_q^n$ with $r=\wt(\bm{x})\le 2$. We show that
		$B_{\bm{x},j}(w)$ depends only on $r$ and $j$.
		
		For a codeword $\bm{c}\in\Aw_w(\C)$, partition the coordinates in
		$\supp(\bm{x})$ into three disjoint sets:
		\begin{align*}
			A&=\{i:c_i=x_i\},\\
			B&=\{i:c_i\ne0,\ c_i\ne x_i\},\\
			C&=\{i:c_i=0\},
		\end{align*}
		where the coordinates are restricted to $\supp(\bm{x})$. Write $a=|A|,\ b=|B|,\ c=|C|$.
		Then $a+b+c=r$.
		Since $\bm{c}$ has weight $w$, we have $d(\bm{x},\bm{c})=w-a+c$.
		By Lemma~\ref{lem:local-pattern}, since $\Aw_w(\C)$ is a
		$q$-ary $2$-design and $a+b+c=r\le2$, the number of codewords
		$\bm{c}\in\Aw_w(\C)$ realizing the above local pattern is a constant
		$\lambda_{a,b,c}$, independent of
		the particular support and nonzero entries of $\bm{x}$.
		
		For fixed $a,b,c$, there are $\frac{r!}{a!b!c!}$
		ways to partition the $r$ coordinates of $\supp(\bm{x})$ into the
		three classes. Therefore,
		\[
		B_{\bm{x},j}(w)=\sum_{\substack{a,b,c\ge0\\a+b+c=r\\w-a+c=j}}
		\frac{r!}{a!b!c!}\lambda_{a,b,c}.
		\]
		The right-hand side depends only on $r=\wt(\bm{x})$ and $j$, and
		not on the particular vector $\bm{x}$.
		
		Consequently, if $\bm{x},\bm{y}\in\F_q^n$ satisfy
		$\wt(\bm{x})=\wt(\bm{y})\le2$, then
		\[
		B_{\bm{x},j}(w)=B_{\bm{y},j}(w),\qquad 0\le j\le n,
		\]
		which completes the proof.
	\end{proof}
	To pass from these local distance counts to the full distance distribution, we use the following elementary correspondence between distance distributions and coset weight distributions.
	\begin{observation}\label{lem: Bxj=A(x+C)}
		Let $\cal C$ be an $[n,k,d]_q$ linear code. For any $\bm{x}\in\F_{q}^{n}$ and any integer $0\le j\le n$, the number of vectors of weight $j$ in the coset $\bm{x}+\C$ equals  $B_{\bm{x},j}$.
	\end{observation}
	\begin{proof}
		Since $\C$ is linear, $-\C=\C$. Hence
		\begin{align*}
			B_{\bm{x},j} &= |\{\bm{c} \in \C : \wt(\bm{x}-\bm{c})=j \}| \\
			&= |\{\bm{v} \in \bm{x}+\C : \wt(\bm{v})=j \}| \\
			&=A_{j}(\bm{x}+\C).\qedhere
		\end{align*}
	\end{proof}
	For MDS codes, the initial coset weight counts determine the remaining weight distribution. Together with Lemma~\ref{lem:local-distance}, this allows us to propagate the	$q$-ary $2$-design property from minimum weight to all	nonzero weights.
	\begin{proposition}\label{prop:MDS-propagation}
		Let $\C$ be an $[n,n-d+1,d]_q$ MDS code with $d\ge2$. If
		$\Aw_d(\C)$ forms a $q$-ary $2$-design,
		then every nonempty set
		$\Aw_w(\C)$ forms a $q$-ary $2$-design.
	\end{proposition}
	\begin{proof}
		Let $s^\perp$ denote the number of distinct nonzero weights in
		$\C^\perp$. Since $\C^\perp$ is an
		$[n,d-1,n-d+2]_q$ MDS code, its nonzero weights belong to $\{n-d+2,n-d+3,\ldots,n\}$.
		Hence $s^\perp\le d-1$.
		We distinguish two cases.
		
		\medskip
		\noindent\textit{Case 1: $s^\perp\le d-2$.}
		In this case, $d-s^\perp\ge2$.
		By the Standard Criterion
		\cite[Theorem~3.6]{WuZhou2026}, every nonempty set
		$\Aw_w(\C)$ forms a $q$-ary $(d-s^\perp)$-design.
		In particular, every nonempty $\Aw_w(\C)$ forms a
		$q$-ary $2$-design, and the assertion follows.
		
		\medskip
		\noindent\textit{Case 2: $s^\perp=d-1$.}
		Fix a vector $\bm{x}\in\F_q^n$ with $\wt(\bm{x})=2$.
		Since $\Aw_d(\C)$ forms a $q$-ary $2$-design, there exists
		a codeword $\bm{c}\in\Aw_d(\C)$ covering $\bm{x}$.
		Therefore, $d(\bm{x},\bm{c})=d-2$, and hence
		\[
		d(\bm{x},\C)\le d-2=s^\perp-1.
		\]
		Thus the minimum weight (or coset weight) of $\bm{x}+\C$ is less than $s^\perp$.
		
		We first determine the initial part of its weight distribution.
		Let
		\[
		0\le j\le s^\perp-1=d-2.
		\]
		If $\bm{c}\in\C$ satisfies $d(\bm{x},\bm{c})=j$, then by the
		triangle inequality,
		\[
		\wt(\bm{c})\le d(\bm{x},\bm{c})+\wt(\bm{x})\le j+2\le d.
		\]
		Since $\C$ has minimum distance $d$, such a codeword is either
		the zero codeword or has weight exactly $d$. Consequently,
		\[
		B_{\bm{x},j}=\delta_{j,2}+B_{\bm{x},j}(d),\qquad 0\le j\le d-2,
		\]
		where $\delta_{j,2}$ is the Kronecker delta.
		By Lemma~\ref{lem:local-distance},
		$B_{\bm{x},j}(d)$ depends only on $\wt(\bm{x})=2$ and $j$.
		Hence $B_{\bm{x},j}$ is independent of the particular vector $\bm{x}$ of weight $2$ for every $0\le j\le s^\perp-1$.
		
		By Observation~\ref{lem: Bxj=A(x+C)}, $B_{\bm{x},j}$ equals the number of vectors of weight $j$ in the coset $\bm{x}+\C$. Since this coset has weight
		less than $s^\perp$, its complete weight distribution is uniquely
		determined by the numbers of vectors of weights
		$1,\ldots,s^\perp-1$
		\cite[Th.~7.5.2]{HuffmanPless}.
		It follows that $B_{\bm{x},j}$ is independent of the particular vector $\bm{x}$ for
		every $0\le j\le n$.
		
		We now prove the desired assertion by induction on the nonzero
		weights of $\C$. The case $w=d$ holds by assumption.
		Suppose that $w>d$ is a nonzero weight of $\C$ and that every
		nonempty $\Aw_v(\C)$ with $d\le v<w$ forms a $q$-ary $2$-design.
		Fix again a vector $\bm{x}$ of weight $2$ and consider
		$B_{\bm{x},w-2}$. If a codeword $\bm{c}\in\C$ contributes to this
		quantity, then
		\[
		\wt(\bm{c})\le d(\bm{x},\bm{c})+\wt(\bm{x})=(w-2)+2=w.
		\]
		Therefore,
		\[
		B_{\bm{x},w-2}(w)=B_{\bm{x},w-2}-\delta_{w-2,2}-\sum_{v=d}^{w-1}B_{\bm{x},w-2}(v).
		\]
		The first term on the right-hand side is independent of $\bm{x}$
		by the coset weight-distribution argument above. For every
		$v<w$ with $\Aw_v(\C)\ne\varnothing$, the induction hypothesis
		and Lemma~\ref{lem:local-distance} imply that
		$B_{\bm{x},w-2}(v)$ is also independent of $\bm{x}$.
		Thus $B_{\bm{x},w-2}(w)$ is constant over all vectors
		$\bm{x}$ of weight $2$.
		
		Finally, for $\bm{c}\in\Aw_w(\C)$ and $\wt(\bm{x})=2$,
		\[
		d(\bm{x},\bm{c})=w-2
		\]
		if and only if $\bm{c}$ covers $\bm{x}$. Hence
		$B_{\bm{x},w-2}(w)$ is exactly the number of codewords in
		$\Aw_w(\C)$ covering $\bm{x}$.
		Since this number is independent of the choice of $\bm{x}$,
		$\Aw_w(\C)$ forms a $q$-ary $2$-design.
		
		By induction, every nonempty $\Aw_w(\C)$ forms a
		$q$-ary $2$-design.
	\end{proof}

	\section{Complete Characterization for DRS Codes}
	
	We first state the multiplicative counting fact that provides the
	necessary condition. Put $G=\F_q^*$ and $n=q-1$. For $1\le r\le n$ and
	$c\in G$, let $N_G(r,c)$ denote the number of $r$-subsets of $G$ whose elements have product $c$.
	
	\begin{lemma}\label{lem:subset-product}
		For $1\le r\le n$, the number $N_G(r,c)$ is independent of $c\in G$ if
		and only if $\gcd(r,n)=1$.  In that case,
		\begin{equation}\label{eq:subset-uniform}
			N_G(r,c)=\frac{1}{n}\binom{n}{r}\qquad(c\in G).
		\end{equation}
	\end{lemma}
	
	\begin{proof}
		If $\gcd(r,n)=1$, the map $h\mapsto h^r$ permutes the cyclic group $G$.
		Given $c_1,c_2\in G$, choose $h$ with $h^r=c_2c_1^{-1}$. Then
		$S\mapsto hS$ is a bijection between $r$-subsets with product $c_1$ and those with product
		$c_2$.
		Hence $N_G(r,c)$ is independent of $c\in G$. Since all
		$\binom{n}{r}$ $r$-subsets are partitioned according to their
		products in $G$, and $|G|=n$, we obtain $N_G(r,c)=\frac{1}{n}\binom{n}{r}$ for every $c\in G$.
		
		Conversely, suppose that $N_G(r,c)$ is independent of $c$, and assume, to the contrary, that $d=\gcd(r,n)>1$. Since $G$ is cyclic, there exists a multiplicative character
		$\chi:G\to\mathbb{C}^{*}$ of order $d$.
		Define
		\[
		A_\chi=\sum_{c\in G}N_G(r,c)\chi(c).
		\]
		Since $N_G(r,c)$ is independent of $c$,  Lemma~\ref{mul-char-sum} gives
		\begin{equation}\label{eq:Achi-zero}
			A_\chi=0.
		\end{equation}
		On the other hand,
		\begin{align}
			A_\chi=\sum_{\substack{S\subseteq G\\|S|=r}}\chi\!\left(\prod_{x\in S}x\right)
			=[x^r]\prod_{a\in G}(1+\chi(a)x), \label{eq:Achi-gen}
		\end{align}
		where $[x^{r}]F(x)$ denotes the coefficient of $x^{r}$ in $F(x)$.
		Since each element of $U_d$ occurs exactly $n/d$ times among the values
		$\chi(a)$ as $a$ ranges over $G$, we have
		\[
		\prod_{a\in G}\bigl(1+\chi(a)x\bigr)=\prod_{\xi^d=1}(1+\xi x)^{n/d}.
		\]
		Using $\prod_{\xi^d=1}(t-\xi)=t^d-1$ and substituting $t=-x^{-1}$, we obtain $\prod_{\xi^d=1}(1+\xi x)=1-(-x)^d$. Consequently,
		\[
		\prod_{a\in G}\bigl(1+\chi(a)x\bigr)=\bigl(1-(-x)^d\bigr)^{n/d}.
		\]
		Since $d\mid r$, write $r=dm$. Then, by \eqref{eq:Achi-gen},
		\begin{align*}
			A_\chi=[x^{dm}]\bigl(1-(-x)^d\bigr)^{n/d}=(-1)^{r+r/d}\binom{n/d}{r/d}.
		\end{align*}
		As $1\le r\le n$, the binomial coefficient is nonzero, and hence
		$A_\chi\neq0$. This contradicts \eqref{eq:Achi-zero}. Therefore,
		$\gcd(r,n)=1$, completing the proof.
	\end{proof}

	The preceding lemma identifies the precise arithmetic condition for uniform subset-product counts.
	In DRS codes, the number of minimum-weight codewords covering a vector of weight $2$ reduces to exactly such a subset-product count, with the prescribed product determined by the ratio of its two nonzero entries. We can therefore translate Lemma~\ref{lem:subset-product} directly into a complete characterization.

	\begin{theorem}\label{thm:min-iff}
		Let $\C_{k}$ be the DRS code defined in \eqref{eq:DRS},  with parameters $[q+1,k,d=q-k+2]_q$, where $q$ is a prime power and $2\le k\le q$.		
		Then $\Aw_d(\C_{k})$ forms a $q$-ary $2$-$\left(q+1,q-k+2,\lambda\right)$ design for some integer $\lambda$ if and only if $\gcd(k-1,q-1)=1$.	Moreover, under this condition, $\lambda=\frac{1}{q-1}\binom{q-1}{k-1}$.
	\end{theorem}
	\begin{proof}
		Suppose first that $\gcd(k-1,q-1)=1$. Then the conclusion follows
		directly from \cite[Theorem~5.5]{WuZhou2026}.
		
		Conversely, for the DRS code $\C_{k}$, with coordinate positions indexed by $\F_{q}\cup\{\infty\}$,	suppose that $\Aw_d(\C_k)$ forms a $q$-ary $2$-design.
		For any $a,b\in\F_q^*$, let $\bm{v}_{a,b}$ be the vector whose entries at positions  $0$ and $1$ are $a$ and $b$, respectively, and whose
		remaining entries are zero. Denote by $\lambda_{a,b}$ the number of
		codewords in $\Aw_d(\C_k)$ that cover $\bm{v}_{a,b}$.
		Since $\Aw_d(\C_k)$ is assumed to form a $q$-ary $2$-design,
		$\lambda_{a,b}$ is independent of the particular choices of $a,b\in\F_q^*$.
		
		Obviously, a minimum-weight codeword of $\C_k$ is determined by a polynomial
		$f\in\mathcal{P}_k[x]$ having exactly $k-1$ zeros in $\F_q\cup\{\infty\}$.

		Separating the cases $\deg(f)=k-1$ and $\deg(f)=k-2$, the counting
		argument in \cite[Theorem~5.5]{WuZhou2026} gives
		\[
		\lambda_{a,b}=N_G(k-1,z),
		\]
		where $z=(-1)^{k-1}f(0)/f(1)=(-1)^{k-1}a/b$.
		As $a$ and $b$ range over $\F_q^*$, the element
		$(-1)^{k-1}a/b$ ranges over all of $\F_q^*$. Hence the independence
		of $\lambda_{a,b}$ from $a$ and $b$ implies that $N_G(k-1,z)$
		is independent of $z\in\F_q^*$.
		By Lemma~\ref{lem:subset-product}, this is possible if and only if $\gcd(k-1,q-1)=1$.
		This proves the necessity.
		
		Finally, under this condition, Lemma~\ref{lem:subset-product} gives
		\[
		N_G(k-1,z)=\frac{1}{q-1}\binom{q-1}{k-1},\qquad z\in\F_q^*,
		\]
		which yields the stated value of $\lambda$.
	\end{proof}

	Combining Theorem~\ref{thm:min-iff} with
	Proposition~\ref{prop:MDS-propagation} yields the main result.
	
	\begin{theorem}\label{thm:all-iff}
		Let $\C_k$ be a DRS code with parameters $[q+1,k,q-k+2]_q$, where
		$2\le k\le q$. The following statements are equivalent:
		\begin{enumerate}
			\item [(1)] $\gcd(k-1,q-1)=1$;
			\item [(2)] $\Aw_{q-k+2}(\C_k)$ forms a $q$-ary $2$-design;
			\item [(3)] every nonempty set $\Aw_w(\C_k)$ forms a $q$-ary
			$2$-design.
		\end{enumerate}
	\end{theorem}
	
	\begin{proof}
		The equivalence of (1) and (2) follows from
		Theorem~\ref{thm:min-iff}. If (1) holds, then
		Proposition~\ref{prop:MDS-propagation} applied to (2) yields (3).
		Finally, (3) immediately implies (2).
	\end{proof}
	
	Theorem~\ref{thm:all-iff} determines exactly when the $q$-ary	$2$-design property holds at every nonzero weight.
	We next make this all-weight characterization quantitative by	determining the corresponding index at each weight.

	\begin{corollary}\label{cor:all-weight-parameters}
		Let $\C_k$ be a DRS code with parameters $[q+1,k,d]_q$, where
		$2\le k\le q$, and suppose that $\gcd(k-1,q-1)=1$.
		Then, for every nonzero weight $w$ occurring in $\C_k$, the set $\Aw_w(\C_k)$ forms a
		$q$-ary $2$-$(q+1,w,\lambda_w)$ design, where
		\[
		\lambda_w=\frac{1}{q-1}\binom{q-1}{w-2}\sum_{j=0}^{w-d}(-1)^j\binom{w-1}{j}q^{\,w-d-j}.
		\]
	\end{corollary}
	
	\begin{proof}
		By Theorem~\ref{thm:all-iff}, $\Aw_w(\C_k)$ forms a $q$-ary $2$-design for every nonzero weight $w$.
		Since $\C_k$ is an MDS code with parameters $[q+1,k,d]_q$, its weight distribution is given by \cite{HuffmanPless}
		\[
		A_w=(q-1)\binom{q+1}{w}\sum_{j=0}^{w-d}(-1)^j\binom{w-1}{j}q^{\,w-d-j}.
		\]
		Applying \eqref{eq:lambda-double-count} with $t=2$ to a $q$-ary $2$-$(q+1,w,\lambda_w)$ design yields
		\[
		\lambda_w=\frac{A_w\binom{w}{2}}{(q-1)^2\binom{q+1}{2}}.
		\]
		Substituting the expression for $A_w$ and using
		\[
		\frac{\binom{q+1}{w}\binom{w}{2}}{\binom{q+1}{2}}=\binom{q-1}{w-2},
		\]
		we obtain the asserted formula.
	\end{proof}

	We conclude this section with two further consequences and extensions	of the above characterization, concerning duality and the possible occurrence of higher-strength $q$-ary designs.

	\begin{remark}
		For $2\le k\le q-1$, the criterion in Theorem~\ref{thm:all-iff} is invariant under duality for DRS codes. Indeed,
		$\C_k^\perp = \C_{q+1-k}$, where $2\le q+1-k\le q-1$, and
		$\gcd(q-k,q-1)=\gcd(k-1,q-1)$. Hence $\C_{k}$ and its dual satisfy the
		all-weight $q$-ary $2$-design property simultaneously.
	\end{remark}
	\begin{remark}
		Let $q\ge 3$. If $\Aw_d(\C_k)$ forms a $q$-ary $3$-design, then
		$d\ge 3$, and hence $2\le k\le q-1$. Since every $q$-ary
		$3$-design is also a $q$-ary $2$-design,
		Theorem~\ref{thm:min-iff} yields $\gcd(k-1,q-1)=1$.
		Moreover, by \eqref{eq:lambda-double-count} and
		\eqref{eq:Ad-DRS}, its index must be $\lambda=\frac{1}{(q-1)^2}\binom{q-2}{k-1}$.
		Therefore, the following arithmetic conditions are necessary:
		\begin{equation}\label{eq:necessary-3design}
			\gcd(k-1,q-1)=1,
			\qquad
			(q-1)^2 \biggm| \binom{q-2}{k-1}.
		\end{equation}
		A direct computation of \eqref{eq:necessary-3design} for all prime
		powers $3\le q<1849$ yields no admissible parameters.
		The smallest prime power not excluded by these necessary conditions is $q=1849=43^2$,
		for which the two dual dimensions $k=828$ and $k=1022$ satisfy
		\eqref{eq:necessary-3design}.
		These conditions are only necessary, however, and do not imply
		that the corresponding minimum-weight codewords form a
		$q$-ary $3$-design.
		Through a series of careful analyses, we
		conjecture that no DRS code with
		$2\le k\le q-1$ holds a $q$-ary $3$-design at minimum weight.
		More generally, our previous investigation \cite{WuZhou2026} found $q$-ary $3$-designs arising from linear codes only in a few
		small cases.
	\end{remark}

	\section{Conclusion}
	
	We have obtained a complete characterization of the $q$-ary	$2$-design property for doubly-extended Reed--Solomon codes. The key step is an exact subset-product criterion over $\F_q^*$,
	which shows that the covering numbers of the minimum-weight codewords are uniform precisely when $\gcd(k-1,q-1)=1$.
	Together with a propagation result for MDS codes, this yields that the same coprimality condition characterizes the $q$-ary $2$-design property at every nonzero weight, with the corresponding indices
	given explicitly by the MDS weight distribution.
	These results resolve the necessity question left open by our	previous work and establish the corresponding all-weight characterization.

	\section*{Acknowledgment}
	This work was supported by the National Natural Science Foundation of China under
	Grant 12571346. OpenAI ChatGPT was used for language and presentation assistance
	in selected passages of Sections~III and~IV.
	All mathematical content and references were independently verified and
	finalized by the authors.

\end{document}